\documentclass[12pt]{article}
\usepackage{hyperref,amsfonts,amssymb,amsmath}
\usepackage{mathtools}

\usepackage{color}

\def\tr{\mathrm{tr}}

\usepackage{soul}

\def\C{\mathbb{C}}
\def\N{\mathbb{N}}
\def\Z{\mathbb{Z}}
\def\Q{\mathbb{Q}}

\def\R{\mathbb{R}}
\def\P{\mathbb{P}}
\def\A{\mathbb{A}}

\def\bq{ \begin{equation} }
\def\eq{ \end{equation} }
\def\ben{ \begin{eqnarray} }
\def\en{ \end{eqnarray} }

\def\frac#1#2{{#1\over #2}}

\def\on#1#2{\mathop{\vbox{\ialign{##\crcr\noalign{\kern2pt}
$\scriptstyle{#2}$\crcr\noalign{\kern2pt\nointerlineskip}
\kern-2pt$\hfil\displaystyle{#1}\hfil$\crcr}}}\limits}

\begin{document}

\title{p-Determinants and monodromy of differential operators}
\author{Maxim Kontsevich,  Alexander Odesskii}
   \date{}
\vspace{-20mm}
   \maketitle
\vspace{-7mm}
\begin{center}
IHES, 35 route de Chartres, Bures-sur-Yvette, F-91440,
France  \\[1ex]
and \\[1ex]
Brock University, 1812 Sir Isaac Brock Way, St. Catharines, ON, L2S 3A1 Canada\\[1ex]
e-mails: \\
\texttt{maxim@ihes.fr}\\
\texttt{aodesski@brocku.ca}
\end{center}

\medskip

\begin{abstract}

We prove that $p$-determinants of a certain class of differential operators can be lifted to power series over $\Q$. We compute these power series in terms of monodromy of the corresponding differential operators.

\end{abstract}

\newpage

\tableofcontents

\newpage

 \section{Introduction}

Consider a differential operator
$$D=\frac{d}{dx}x(x-1)(x-t)\frac{d}{dx}+x=x(x-1)(x-t)\frac{d^2}{dx^2}+(3x^2-2(t+1)x+t)\frac{d}{d x}+x$$ 
where $t$ is a parameter. This operator has regular singular points at $x=0,1,t,\infty$ with a monodromy matrix $\left( \begin{array}{cc} 1 & 1\\ 0 & 1\end{array}  \right)$ at each singular point. Therefore, the monodromy group of the equation  $Df(x)=0$ is a 
subgroup of $SL_2(\C)$. 

Let us discuss the eigenvalues of the  monodromy matrix   for a loop with $0,t$ inside and $1,\infty$ outside  for $ |t|\ll 1$. Denote these eigenvalues by 
\begin{equation}\label{lam}
\exp(\pm2\pi i \lambda). 
\end{equation}
It is clear that an eigenfunction $g_{\pm}(x)$ corresponding to the eigenvalue $\exp(\pm 2\pi i\lambda)$ can be written as

\begin{equation}\label{eig}
g_{\pm}(x)=x^{\pm\lambda}\sum_{l\in\Z}r^{\pm}_lx^l.
\end{equation}

Here we assume that $|t|<|x|<1$ and $\lambda, r^{\pm}_l$ are certain analytic functions in $t$ for $0\le |t|\ll 1$. Moreover, in the limit $t=0$ the differential equation $D f(x)=0$ has two solutions $f(x)=\frac{1}{x}$ and $f(x)=\frac{\ln(1-x)}{x}$ with $\lambda=0$. Therefore\footnote{In principle,  we could have a function $\lambda(t)$  depending analytically on $\sqrt{t}$. However, one can show by direct calculations that $\lambda(t)$ is analytic in $t$, see a formula (\ref{lambda}) below and Example 1 in the last Section.}, for small $t$ the functions $\lambda, r^{\pm}_l$  have power series expansion in $t$. 

Consider a reduction of the operator $D$ modulo a prime $p$.  Let $V$ be a free module with a basis $1,x,x^2,...$ over  the ring $\Z/p\Z[t]$. Let $V_p\subset V$ be a  submodule of $V$ with a basis $x^p,x^{p+1},...$. It is clear that $V_p$ is invariant with respect to $D$ modulo $p$. Denote by $D_p$ the corresponding linear operator on $V/V_p$. In the basis $1,x,...,x^{p-1}$ of $V/V_p$ we have 
\begin{equation}\label{pmat}
D_p=\left( \begin{array}{cccccccccc}
0 &  1^2t & & 0\\
 1^2 &  1\cdot2(t+1) &  2^2t &  &\\
&  2^2 &  2\cdot3(t+1) &   3^2t \\
 &  & \ddots &  \ddots  &   \\
 &&&&&\\
 0&&  (p-1)^2& (p-1)p(t+1)\end{array}  \right) 
\end{equation}
Denote $Det_p(D)=\det D_p$, it is called a $p$-determinant\footnote{See \cite{K} for general definition and discussion of $p$-determinants.} of the differential operator $D$. It is clear that $Det_p(D)\in\Z/p\Z[t]$.

{\bf Claim.} There exists a power series $h(t)=\sum_{n=0}^{\infty}c_nt^n\in\Q[[t]]$ independent of $p$ such that for each prime $p$ we have\footnote{This means that  $c_0,c_1,\dots,c_{p-1}$ are rational numbers without $p$ in denominator.}
 $c_0,c_1,\dots,c_{p-1}\in \Z_p\cap \Q$ and 
$$Det_p(D)=\sum_{i=0}^{p-2} c_i({\rm mod}\,\, p) \cdot t^i-t^{p-1}.$$
Informally, we can write
$$Det_p(D)=h(t)~\text{mod}~p,~t^{p-1}.$$ Moreover, $h(t)$ has  a  nonzero radius of convergence $r$ and if $|t|<r$, then we have 
$$h(t)=\lambda^2$$
where $\lambda$ is defined by (\ref{lam}).

 First few coefficients of the series $h(t)$ are
$$h(t)=\frac{1}{4}t^2+\frac{1}{24}t^3+\frac{101}{576}t^4+\frac{239}{17280} t^5+\frac{19153}{115200}t^6-\frac{1516283}{72576000}t^7+\frac{ 23167560743}{121927680000}t^8 +\dots.$$
Correspondingly, we have
\begin{equation}
 \label{lambda}  
\lambda(t)=\frac{1}{2}t +\frac{1}{24}t^2+\frac{25}{144}t^3 -\frac{11}{17280}t^4 + \frac{70591}{518400}t^5 - \frac{774601}{24192000}t^6 + \frac{2215989011}{15240960000}t^7+ \dots
\end{equation}
Both series $h(t),\lambda(t)$ have quite unusual arithmetic properties, e.g. seem to have \emph{zero} radius of convergence $p$-adically for all primes  $p>2$. Also, it is computationally very hard to calculate rational numbers $c_i$. The last known to us value $c_{251}$ is the ratio of two integers each  of size $\sim 10^{12000}$.

In this paper we consider $p$-determinants of a class of differential operators\footnote{Our previous concrete operator $D$ can be reduced to this canonical form after a multiplication by a function of $t$. See Example 1 in the last Section for details.} of a form 
\begin{equation}
  D=(x\partial_x-l_1)...(x\partial_x-l_n)+\sum_{\substack{0\leq j\leq m,\\ 0\leq i}}t_{i,j}x^i\partial_x^j,~~~~~~~\partial_x:=\frac{d}{dx} \label{family}  
\end{equation}
for some $m\in\N,~l_1,...,l_n\in\Z$ and $l_1\leq l_2\leq...\leq l_n$. Here $t_{i,j}$ are parameters which considered as small  for $i\leq j$. Alternatively
 one can interpret $D$ as a small perturbation of a differential operator with the symbol $x^n\frac{d^n}{dx^n}$,  a  regular singularity at $x=0$, and with the  unipotent monodromy around $x=0$.  We prove that the $p$-determinant of $D$ can be lifted to an  independent of $p$ power series in small parameters. Moreover, we prove that if $m\leq n$, then
$$Det_p(D)=\det \frac{\log M_D}{2\pi i}~ \text{~mod~}~ p \text{~~and large powers of small parameters}$$
where $M_D$ is a monodromy matrix of $D$ with respect to a small fixed loop around zero.  

We formulate our results in terms of a commutative ring $R$ and a proper ideal $I\subset R$. Intuitively, $I$ consists of small elements of $R$. For example, $R=\Z[t]$ and $I=(t)$. We use the notation $\hat{R}_I$ for $I$-adic completion of $R$ 
which is a generalization of the ring of power series $\Z[[t]]$. We also use reduction of $R$ modulo a prime $p$ which is $R/pR$. Our lifting of a $p$-determinant is an element of $\hat{R}_I\widehat{\otimes}_{\Z}\Q$. In our example this is just the ring  $\Q[[t]]$ of power series in $t$ with rational coefficients.

The proofs of our main results are based on certain combinatorial lemmas  concerning  infinite matrices. 
These lemmas are collected in Section 2. Main result of this Section is Lemma 3, formula (\ref{limp4}).

Section 3 is devoted to $p$-determinants of certain class of infinite matrices with polynomial coefficients. 
Main result of this Section is the formula (\ref{lift1}).

In Section 4 we prove an algebraic version of our main result about $p$-determinants of differential operators (see formula (\ref{pdetinf})), and in Section 5 we prove its analytic version. Section 6 is devoted to examples. We discuss further our main example from this Introduction and another example related to Hasse invariants of elliptic curves studied in \cite{M}.

 \section{Matrix lemmas}
 
Let $R$ be a commutative ring  and $I\subset R$ be a proper ideal. 
Throughout the paper we will make the following assumption on the quotient ring $R_0:=R/I$: both maps 
\begin{equation}
    R_0\to R_0\otimes_\Z \Q,~~~~R_0\to \left(\prod_{p=2,3,5,\dots} R_0/pR_0\right) /\bigoplus_{p=2,3,5,\dots} R_0/pR_0\quad \label{ring_assumption}
\end{equation}
are inclusions. In fact, the first assertion follows from the second one. 
This assumption holds e.g. for any finitely-generated ring which is torsion-free as a $\Z$-module.

We denote by $\hat{R}_I$ the $I$-adic completion\footnote{For example $R=\Z[t]$ and $I=(t)$. In this case $\hat{R}_I=\Z[[t]]$.} of $R$.
 
Let $A=(a_{ij}),~i,j\in\Z$ be an infinite matrix where $a_{ij}\in R$. In this Section we assume that 
\begin{itemize}
\item
  $A$ has only finitely many nonzero diagonals above the main diagonal:  $a_{ij}=0$ if $j-i>m$ for some fixed $m\in\N$.
\item
  $a_{ii}=1$ modulo $I$ for all $i\in\Z$.
\item
  $a_{ij}\in I$ if $i<j$.
\end{itemize}

It will be convenient to introduce notation 
$$\overline{a}_{ij}:=\begin{cases}a_{ij} &\text{ if }i\ne j\\ a_{ii}-1&\text{ if } i=j\end{cases}$$
hence $\overline{a}_{ij}\in I$ for $i\le j$.

For any integers $a\leq b$ let $A_{ab}$ be a $b-a+1 \times b-a+1$ submatrix of $A$ defined by $A_{ab}=(a_{ij}),~a\leq i,j\leq b$. 

{\bf Lemma 1.} $\log \det A_{ab}$ is a well-defined element of the ring $\hat{R}_I\hat{\otimes}_{\Z}\Q$. 

{\bf Proof.} We have $$\log \det A_{ab}=\tr\log A_{ab}=\sum_{k>0} \frac{(-1)^{k-1}}{k} \tr(A_{ab}-1)^k$$ and we obtain
\begin{equation}\label{ln}
\log \det A_{ab}=\sum_{\substack{a\leq i_1,...,i_k\leq b,\\i_2-i_1,i_3-i_2,...,i_1-i_k\leq m}}\frac{(-1)^{k-1}}{k}\overline{a}_{i_1,i_2}\overline{a}_{i_2,i_3}...\overline{a}_{i_k,i_1}.
\end{equation}
Notice that for each $k$ there are only finitely many terms of the form $\overline{a}_{i_1,i_2}\overline{a}_{i_2,i_3}...\overline{a}_{i_k,i_1}$ in the r.h.s. of (\ref{ln}), this number is bounded from above by $(b-a+1)^k$ because 
$a\leq i_1,...,i_k\leq b$. It suffice to prove that $\overline{a}_{i_1,i_2}\overline{a}_{i_2,i_3}...\overline{a}_{i_k,i_1}\in I^{\lceil\frac{k}{m+1}\rceil}$. This is done in the proof  of the next lemma. $\square$

{\bf Lemma 2.} Let $a,b$ be integers such that $a+b>0$. The $I$-adic limit 

\begin{equation}\label{lim}
\lim_{N,M\to +\infty} \log\Big(\frac{\det A_{-M,-b}\det A_{a,N}}{\det A_{-M,N}}\Big)
\end{equation}
is a well defined element of the ring $\hat{R}_I\hat{\otimes}_{\Z}\Q$.  
 
{\bf Proof.} Using (\ref{ln}) and our assumptions about $A$ we can write this limit as follows
\begin{equation}\label{lim1}
\lim_{N,M\to +\infty} \log\Big(\frac{\det A_{-M,-b}\det A_{a,N}}{\det A_{-M,N}}\Big)=\sum_{\substack{i_1,...,i_k\in\Z,\\i_2-i_1,i_3-i_2,...,i_1-i_k\leq m,\\ \min(i_1,...,i_k)<a,~\max(i_1,...,i_k)>-b}}\frac{(-1)^k}{k}\overline{a}_{i_1,i_2}\overline{a}_{i_2,i_3}...\overline{a}_{i_k,i_1}.
\end{equation}
Notice that for each $k$ there are only finitely many terms of the form $\overline{a}_{i_1,i_2}\overline{a}_{i_2,i_3}...\overline{a}_{i_k,i_1}$ in the r.h.s. of (\ref{lim1}). Indeed, in addition to the inequalities 
$$\min(i_1,...,i_k)<a,~\max(i_1,...,i_k)>-b$$
we also have $$\max(i_1,...,i_k)-\min(i_1,...,i_k)<km$$
because $i_2-i_1,i_3-i_2,...,i_k-i_{k-1},i_1-i_k\leq m$. It follows from these inequalities that $$-b-km<i_1,...,i_k<a+km$$
and, therefore, the number of monomials is bounded from above by $(a+b+2km)^k$. 

It suffice to prove that $\overline{a}_{i_1,i_2}\overline{a}_{i_2,i_3}...\overline{a}_{i_k,i_1}\in I^{\lceil\frac{k}{m+1}\rceil}$. 

Indeed, we have a sum of $k$ integers equal to zero
$$(i_2-i_1)+(i_3-i_2)+\dots+(i_1-i_k)=0 $$
and each summand belongs to $\Z_{\le m}$. Denote by $N$ the number of summands which belong to $\{0,\dots,m\}$, the remaining $k-N$ summands belong to $\Z_{\le -1}$. Hence we have
$$ 0\le Nm+(k-N)\cdot(-1)\quad \iff N\ge \frac{k}{m+1}.\quad\square$$ 
 
{\bf Lemma 3.} Let  $a+b>0$ and $p$ is another positive integer such that $a,b \ll p$. Then 
\begin{multline}\label{limp}
\lim_{N,M\to +\infty} \log\Big(\frac{\det A_{-M,N}}{\det A_{-M,-b}\det A_{a,N}}\Big)+\log \det A_{a,p-b}=\\
\sum_{\substack{ 1\le k<p,\\ \max(i_1,...,i_k)=-b,\\i_2-i_1,...,i_1-i_k\leq m}}
  \frac{(-1)^{k-1}}{k}\sum_{\gamma=1} \overline{a}_{i_1+\gamma,i_2+\gamma}\overline{a}_{i_2+\gamma,i_3+\gamma}...\overline{a}_{i_k+\gamma,i_1+\gamma} \text{ modulo } I^{\lceil\frac{p+1-a-b}{m(m+1)}\rceil}
\end{multline} 
Notice that this formula can be rewritten as 
\begin{equation}\label{limp4}
\det A_{a,p-b}=\lim_{N,M\to +\infty} \frac{\det A_{-M,-b}\det A_{a,N}}{\det A_{-M,N}}\times 
\end{equation}
$$ \times \exp\sum_{\substack{1\leq\gamma\leq p\\ k\ge 1,\\ \max(i_1,...,i_k)=-b,\\i_2-i_1,...,i_1-i_k\leq m}}\frac{(-1)^{k-1}}{k}\overline{a}_{i_1+\gamma,i_2+\gamma}\overline{a}_{i_2+\gamma,i_3+\gamma}...\overline{a}_{i_k+\gamma,i_1+\gamma} \times \exp\phi,~~~\phi\in I^{\lceil\frac{p+1-a-b}{m(m+1)}\rceil} $$

{\bf Proof.}  We will use the following observation: any nonempty sequence of integers can be uniquely written as 
$$i_1+\gamma,i_2+\gamma,\dots, i_k+\gamma$$
where $\max(i_1,\dots,i_j)=-b$ and $\gamma\in \Z$.

For each sequence $i_1,\dots,i_k\in \Z$ with $k<p$ such that $$\max(i_1,...,i_k)=-b,\quad i_2-i_1,...,i_1-i_k\leq m$$ the terms  
$$\frac{(-1)^{k-1}}{k} \overline{a}_{i_1+\gamma,i_2+\gamma}\overline{a}_{i_2+\gamma,i_3+\gamma}...\overline{a}_{i_k+\gamma,i_1+\gamma}$$
appear in the expansion of $\lim_{N,M\to +\infty}\log\Big(\frac{\det A_{-M,N}}{\det A_{-M,-b}\det A_{a,N}}\Big) $
only when
$$\mathbf{min}+\gamma<a, \quad \mathbf{max}+\gamma>-b$$
where 
$$\mathbf{min}:=\min(i_1,\dots,i_k),\quad \mathbf{max}:=\max(i_1,\dots,i_k)=-b$$
Similarly, the same terms appear in the expansion of $\log \det A_{a,p-b}$ only when
$$\mathbf{min}+\gamma\ge a, \quad \mathbf{max}+\gamma\le p-b$$
The union of sets
$$\{\gamma\in \Z|\,\mathbf{min}+\gamma<a, \quad \mathbf{max}+\gamma>-b\}\cup \{\gamma\in \Z|\,\mathbf{min}+\gamma\ge a, \quad \mathbf{max}+\gamma\le p-b\}   $$
is disjoint and coincides with $\{1,\dots,p\}$ if $a-\mathbf{min}\le p+1$.

The latter condition can be rewritten as $(a+b)+(\mathbf{max}-\mathbf{min})\le (p+1)$. If it does not hold, then $p+1-a-b<\mathbf{max}-\mathbf{min}$. On the other hand, $\mathbf{max}-\mathbf{min}<km$ because $i_2-i_1,i_3-i_2,...,i_k-i_{k-1},i_1-i_k\leq m$. It follows from these inequalities that $k>\frac{p+1-a-b}{m}$ 
and therefore $ \overline{a}_{i_1+\gamma,i_2+\gamma}\overline{a}_{i_2+\gamma,i_3+\gamma}...\overline{a}_{i_k+\gamma,i_1+\gamma}\in I^{\lceil\frac{p+1-a-b}{m(m+1)}\rceil}$.  $\square$

\section{Matrix p-determinants depending on $\varepsilon$}

Let $B(\varepsilon)=(b_{ij}(\varepsilon)),~i,j\in\Z$ be an infinite matrix where\footnote{Here $R[\varepsilon]$ is a ring of polynomials in an independent variable $\varepsilon$  with coefficients in $R$. Later variable $\varepsilon$ will appear as a parameter for the twist (conjugation) by $x^\varepsilon$ of a differential operator.} $b_{ij}(\varepsilon)\in R[\varepsilon]$. In this Section we assume that 
\begin{itemize}
\item
  $B(\varepsilon)$ has only finitely many nonzero diagonals above the main diagonal:  $b_{ij}(\varepsilon)=0$ if $j-i>m$ for some fixed $m\in\N$.
\item
  $b_{ii}(\varepsilon)=(i+\varepsilon-l_1)(i+\varepsilon-l_2)...(i+\varepsilon-l_n)$ modulo $I[\varepsilon]$ for all $i\in\Z$. Here $l_1,l_2,...,l_n\in\Z$ and $l_1\leq l_2\leq...\leq l_n$. 
\item
  $b_{ij}(\varepsilon)\in R[\varepsilon]$ and, moreover $b_{ij}(\varepsilon)=q_{j-i}(\varepsilon+j)$ for $i,j\in\Z$ such that $j-i\leq m$. Here $q_{i}(\varepsilon)\in R[\varepsilon]$.
\item $q_{i}(\varepsilon)\in I[\varepsilon]$ if $i>0$.
\end{itemize}

For any integers $a\leq b$ let $B_{ab}(\varepsilon)$ be a $b-a+1 \times b-a+1$ submatrix of $B(\varepsilon)$ defined by $B_{ab}(\varepsilon)=(b_{ij}(\varepsilon)),~a\leq i,j\leq b$.

In the sequel of this Section we assume that  $a+b>0$ and $p$ is an arbitrary prime number such that $a,b,n\ll p$. We are interested in reductions of $\det B_{a,p-b}(\varepsilon)$ modulo $p$.  

Consider the ring $$\widetilde{R}:=R[\varepsilon,\frac{1}{\varepsilon},\frac{1}{\varepsilon\pm 1},\frac{1}{\varepsilon\pm 2},...]$$
We have 
$$\widetilde{R}\otimes \Z/p\Z=(R/pR)[\varepsilon,\frac{1}{\varepsilon},\frac{1}{\varepsilon- 1},\frac{1}{\varepsilon - 2},\dots,\frac{1}{\varepsilon-(p-1)}]\subset R/pR((\epsilon))$$

Define an $R$-linear operator $$T_{\varepsilon}:\widetilde{R}=R[\varepsilon,\frac{1}{\varepsilon},\frac{1}{\varepsilon\pm 1},\frac{1}{\varepsilon\pm 2},...]\to \frac{1}{\varepsilon}R[\frac{1}{\varepsilon}]\subset \widetilde{R}$$
as follows:
$$T_{\varepsilon}(\varepsilon^i)=0,~i=0,1,2,..., $$$$T_{\varepsilon}\Big(\frac{1}{(\varepsilon-c)^j}\Big)=\frac{1}{\varepsilon^j},~c\in\Z,~j=1,2,...$$

{\bf Lemma 4.} Let $f(\varepsilon)\in R[\varepsilon,\frac{1}{\varepsilon},\frac{1}{\varepsilon\pm 1},\frac{1}{\varepsilon\pm 2},...]$. Then 
$$f(\varepsilon+1)+...+f(\varepsilon+p)= T_{\varepsilon}(f(\varepsilon))+\phi~\text{mod}~p$$
where $\phi\in \varepsilon^{p-l} R/pR[[\varepsilon]]$ and $l$ is the largest multiplicity of roots of denominator of $f(\varepsilon)$.

{\bf Proof.} Using partial fractions representation of $f(\varepsilon)$ we reduce this statement to the identities
$$(\varepsilon-c+1)^l+...+(\varepsilon-c+p)^l=0~~~\text{mod}~p,~~~\text{for~} l=0,1,2,...,$$
$$\Bigg(\frac{1}{(\varepsilon-c+1)^l}+...+\frac{1}{(\varepsilon-c+p)^l}-\frac{1}{\varepsilon^l}\Bigg)~~~\text{mod}~p ~~\in~\varepsilon^{p-l}\Z/p\Z[[\varepsilon]],~~~\text{for~}0<l<p\text{~and}~c\in\Z$$
which are clear. $\square$
 
{\bf Theorem 1.} The following identity holds 
\begin{equation}\label{liftgen}
\det B_{a,p-b}(\varepsilon)=(-\varepsilon)^n\lim_{N,M\to +\infty} \frac{\det B_{-M,-b}(\varepsilon)\det B_{a,N}(\varepsilon)}{\det B_{-M,N}(\varepsilon)}\times
\end{equation} 
$$\times\exp T_{\varepsilon}\left(\lim_{N\to+\infty}\text{tr}(\log A_{-N,-b}(\varepsilon) -\log A_{-N,-b-1}(\varepsilon))\right)\times \exp(\phi),\quad \phi\in pR+I^{\lceil Cp\rceil }\hat{\otimes}\Q+\epsilon^{\lceil Cp\rceil}I\hat{\otimes}\Q$$
where $A(\varepsilon)=(a_{ij}(\varepsilon)),~i,j\in\Z$ is an infinite matrix with entries $$a_{ij}(\varepsilon)=\frac{b_{ij}(\varepsilon)}{(j+\varepsilon-l_1)(j+\varepsilon-l_2)...(j+\varepsilon-l_n)}$$ and the operator $T_{\varepsilon}$ 
is extended to the ring $\hat{R}_I[\varepsilon,\frac{1}{\varepsilon},\frac{1}{\varepsilon\pm 1},\frac{1}{\varepsilon\pm 2},...]\hat{\otimes}_{\Z}\Q$ in the natural way. Here $C\in\R_{>0}$ is independent of $p$.

{\bf Proof.} We have $$\det B_{a,p-b}(\varepsilon)=\prod_{a\leq j\leq p-b} b^0_{jj}(\varepsilon)~\det A_{a,p-b}(\varepsilon)$$
where $$b^0_{jj}=(j+\varepsilon-l_1)(j+\varepsilon-l_2)...(j+\varepsilon-l_n).$$ Notice that $b^0_{jj}=b_{jj}$ mod $I$.
Using the formula (\ref{limp4}) we obtain
$$\det B_{a,p-b}(\varepsilon)=\prod_{a\leq i\leq p-b} b^0_{ii}(\varepsilon)\lim_{N,M\to +\infty} \frac{\det A_{-M,-b}(\varepsilon)\det A_{a,N}(\varepsilon)}{\det A_{-M,N}(\varepsilon)}\times$$
$$\times\exp\sum_{\substack{1\leq\gamma\leq p,\\k>1,\\\max(i_1,...,i_k)=-b,\\i_2-i_1,...,i_1-i_k\leq m}}\frac{(-1)^{k-1}}{k}\overline{a}_{i_1+\gamma,i_2+\gamma}(\varepsilon)\overline{a}_{i_2+\gamma,i_3+\gamma}(\varepsilon)...\overline{a}_{i_k+\gamma,i_1+\gamma}(\varepsilon)\times\exp(\phi_1),\quad \phi_1\in I^{\lceil Cp\rceil}\hat{\otimes}\Q .$$
The first part of this formula can be written  as 
$$\prod_{-b<i\leq p-b} b^0_{ii}(\varepsilon)\lim_{N,M\to +\infty} \frac{\det B_{-M,-b}(\varepsilon)\det B_{a,N}(\varepsilon)}{\det B_{-M,N}(\varepsilon)}$$
and the product of $b^0_{ii}(\varepsilon)$ is equal to $(-\varepsilon)^n$ up to multiplication by $\exp(\phi_2),\quad \phi_2\in (p)+(\epsilon^{p})$. The second part is equal to 
$$\exp T_{\varepsilon}\left( \lim_{N\to+\infty}\text{tr}(\log A_{-N,1}(\varepsilon) -\log A_{-N,0}(\varepsilon))\right)\times\exp(\phi_3),\quad \phi_3\in  pI+\epsilon^{\lceil Cp\rceil}I\hat{\otimes} \Q$$
by Lemma 4.  $\square$

To compute the r.h.s. of the formula (\ref{liftgen}) more explicitly we need the following formal version 
of the Weierstrass preparation theorem\footnote{This statement might be well known to the experts but we were unable to find it in the literature. A very similar but not suitable for us statement can be found in  Bourbaki, Commutative Algebra,  Chapter 7, Section 3.8. The Bourbaki's proof is similar to our Proof 1.}.

{\bf Lemma 5.} Let $$q(\varepsilon)=\sum_{j=0}^{\infty}q_j\varepsilon^j$$
be an element of $\hat{R}_I[[\varepsilon]]$ such that $q_0,q_1,...,q_{n-1}\in I$ and $q_n=1$ modulo $I$. Then there exists a unique 
factorization\footnote{In the sequel we will refer to the polynomial $\varepsilon^n-w_1\varepsilon^{n-1}+...+(-1)^nw_n$ as to the Weierstrass polynomial of $q(\varepsilon)$.}

\begin{equation}\label{f}
q(\varepsilon)=(\varepsilon^n-w_1\varepsilon^{n-1}+...+(-1)^nw_n)q_{\text{inv}}(\varepsilon)
\end{equation}
where $w_1,...,w_n\in I\hat{R}_I$ and $q_{\text{inv}}(\varepsilon)=1$ modulo $I\hat{R}_I+(\varepsilon)$ is an invertible 
element of the ring $\hat{R}_I[[\varepsilon]]$.

We suggest two proofs: the first one is simpler and the second one gives an explicit formula for the Weierstrass polynomial.
\newpage
{\bf Proof 1.}

Let $q_{\text{inv}}(\varepsilon)=v_0+v_1\varepsilon+v_2\varepsilon^2+...$ Equating coefficients at powers of $\varepsilon$ in (\ref{f}) we obtain
\begin{equation}\label{eq1}
\begin{array}{l}
w_nv_0=(-1)^nq_0 \\
w_nv_1-w_{n-1}v_0=(-1)^nq_1\\
............................................\\
w_nv_{n-1}-w_{n-1}v_{n-2}+...+(-1)^{n-1}w_1v_0=(-1)^nq_{n-1}
\end{array}
\end{equation}
by equating coefficients at $1,\varepsilon,...\varepsilon^{n-1}$ and
\begin{equation}\label{eq2}
\begin{array}{l}
v_0=q_n+w_1v_1-w_2v_2+...+(-1)^{n-1}w_nv_n \\
v_1=q_{n+1}+w_1v_2-w_2v_3+...+(-1)^{n-1}w_nv_{n+1}\\
v_3=q_{n+2}+w_1v_3-w_2v_4+...+(-1)^{n-1}w_nv_{n+2}\\
............................................
\end{array}
\end{equation}
by equating coefficients at $\varepsilon^n,\varepsilon^{n+1},\varepsilon^{n+2},...$.

Reducing (\ref{eq2}) modulo $I$ we obtain $v_i=q_{n+i}$ mod $I$ for $i=0,1,2,...$. In particular $v_0=q_n$ modulo $I$ is invertible in $R/I$. Reducing (\ref{eq1}) modulo $I^2$ we obtain a system of linear equations for 
$w_n,w_{n-1},...,w_1$ which has unique solution modulo $I^2$ because $v_0$ is invertible. For example, 
$w_n=(-1)^nq_n^{-1}q_0\in I/I^2$ modulo $I^2$. 
To prove the Lemma it suffice to show that for each $k=1,2,3,...$ there exist unique elements $v_0,v_1,v_2,...\in R/I^k$ and $w_n,...w_1\in I/I^{k+1}$ such that equations (\ref{eq2}) hold modulo $I^k$ and 
equations (\ref{eq1}) hold modulo $I^{k+1}$. This is already done for $k=1$. Using induction by $k$ and assuming
 that we already know $v_0,v_1,v_2,...\in R/I^k$ and $w_n,...w_1\in I/I^{k+1}$ we reduce (\ref{eq2}) by $I^{k+1}$ and find $v_0,v_1,v_2,...\in R/I^{k+1}$. After that we reduce (\ref{eq1}) by $I^{k+2}$ and find 
 $w_n,...w_1\in I/I^{k+2}$. $\square$

{\bf Proof 2.}
Define two $\hat{R}_I$-linear endomorphisms $L_{<0},~L_{\geq 0}$ of the $\hat{R}_I$-module $\hat{R}_I[[\varepsilon^{-1},\varepsilon]]$ by 
$$L_{<0}\Big(\frac{1}{\varepsilon^i}\Big)=\frac{1}{\varepsilon^i},~i=1,2,...\text{~and~} L_{<0}(\varepsilon^j)=0,~j=0,1,...$$
$$L_{\geq 0}\Big(\frac{1}{\varepsilon^i}\Big)=0,~i=1,2,...\text{~and~} L_{\geq 0}(\varepsilon^j)=\varepsilon^j,~j=0,1,...$$
 Define elements $q_{\pm}(\varepsilon)\in \hat{R}_I[[\varepsilon^{-1},\varepsilon]]$ by 
\begin{equation}\label{qpm}
q_{-}(\varepsilon)=\exp \Big(L_{<0} \Big( \log \frac{q(\varepsilon)}{\varepsilon^n}\Big)\Big),~~~~q_{+}(\varepsilon)=\exp \Big(L_{\geq 0} \Big( \log \frac{q(\varepsilon)}{\varepsilon^n}\Big)\Big)
\end{equation} 
where $\log \frac{q(\varepsilon)}{\varepsilon^n}$ is understood as 
$$\log \frac{q(\varepsilon)}{\varepsilon^n}=\sum_{k> 0}\frac{(-1)^{k-1}}{k}\Big(\frac{q(\varepsilon)}{\varepsilon^n}-1\Big)^k$$
$$=\sum_{\substack{i,j\geq 0,\\i+j>0}}\frac{(-1)^{i+j-1}(i+j-1)!}{i!j!}\Big(\frac{q_0}{\varepsilon^n}+\frac{q_1}{\varepsilon^{n-1}}+...+\frac{q_{n-1}}{\varepsilon}\Big)^i(q_n-1+q_{n+1}\varepsilon+q_{n+2}\varepsilon^2...)^j$$
and this is a well defined element of $\hat{R}_I\hat{\otimes}_{\Z}\Q[[\varepsilon^{-1},\varepsilon]]$ because $q_0,...,q_{n-1}\in I$ 
and $q_n=1$ mod $I$.

We have 
$$q_{-}(\varepsilon)q_{+}(\varepsilon)=\exp \Big(L_{<0} \Big( \log \frac{q(\varepsilon)}{\varepsilon^n}\Big)+L_{\geq 0} \Big( \log \frac{q(\varepsilon)}{\varepsilon^n}\Big)\Big)=\frac{q(\varepsilon)}{\varepsilon^n}$$
because $L_{<0}+L_{\geq 0}=1$. Moreover, $q_{+}(\varepsilon)$ is an invertible element of the ring $\hat{R}_I[[\varepsilon]]$ because $q(\varepsilon)=\varepsilon^n+q_{n+1}\varepsilon^{n+1}+...$ modulo $I$ and, therefore $q_{+}(\varepsilon)=1+q_{n+1}\varepsilon+...$ modulo $I$.

We have $q_{-}(\varepsilon)=\frac{q(\varepsilon)}{\varepsilon^nq_{+}(\varepsilon)}$ which shows that 
$q_{-}(\varepsilon)$ does not contain monomials $\frac{1}{\varepsilon^{n+1}},\frac{1}{\varepsilon^{n+2}},...$.
On the other hands, it follows from (\ref{qpm}) that $q_{-}(\varepsilon)$ is a power series in $\frac{1}{\varepsilon}$ with the constant term equal  to one. Therefore, $q_{-}(\varepsilon)$ is a polynomial in
$\frac{1}{\varepsilon}$ of degree less than  or equal to  $n$ and with the constant term equal  to one. This proves the existence of a factorization (\ref{f}) if we set
\begin{equation}\label{qpm1}
\varepsilon^n-w_1\varepsilon^{n-1}+...+(-1)^nw_n=\varepsilon^nq_{-}(\varepsilon) \text{~~~and~~~}  q_{\text{inv}}(\varepsilon)=q_{+}(\varepsilon).
\end{equation}
 To prove uniqueness we apply the operators 
$$f\mapsto \exp(L_{<0}(\log f)),~~~~f\mapsto \exp(L_{\geq 0}(\log f))$$
to the l.h.s and the r.h.s. of the equation
$$\frac{q(\varepsilon)}{\varepsilon^n}=\Big(1-\frac{w_1}{\varepsilon}+...+(-1)^n\frac{w_n}{\varepsilon^n}\Big)q_{\text{inv}}(\varepsilon)$$
and see that if a factorization (\ref{f}) exists, then the Weierstrass polynomial $\varepsilon^n-w_1\varepsilon^{n-1}+...+(-1)^nw_n$ and the invertible power series $q_{\text{inv}}(\varepsilon)$ are given by (\ref{qpm1}) and, therefore, unique.
$\square$

By Lemma 5 we have for large $N,M$ 
$$\det B_{-M,N}(\varepsilon)=(\varepsilon^n-w_{1,-M,N}\varepsilon^{n-1}+...+(-1)^{n}w_{n,-M,N})Q_{-M,N}$$
where  $w_{i,-M,N}\in I  \hat{R}_I\hat{\otimes}_{\Z}\Q$  and $Q_{-M,N}\in\hat{R}_I\hat{\otimes}_{\Z}\Q[[\varepsilon]]$ is  an invertible element in $\hat{R}_I\hat{\otimes}_{\Z}\Q[[\varepsilon]]$. Indeed, 
$$\det B_{-M,N}(\varepsilon)=\prod_{-M\leq i\leq N} b_{ii}(\varepsilon) \text{~modulo~} I=\prod_{-M\leq i\leq N} (i+\varepsilon-l_1)(i+\varepsilon-l_2)...(i+\varepsilon-l_n) \text{~modulo~} I$$
and so $\det B_{-M,N}(\varepsilon)=\varepsilon^n(k_0+k_1\varepsilon+...)$ mod $I$ where $k_0$ is a nonzero integer if $M,N>|l_1|,|l_n|$.

{\bf Theorem 2.} There exists $I$-adic limit
$$w(\varepsilon)=\lim_{N,M\to+\infty}(\varepsilon^n-w_{1,-M,N}\varepsilon^{n-1}+...+(-1)^{n}w_{n,-M,N}).$$
Moreover, we have 
\begin{equation}\label{w}
\varepsilon^n\exp T_{\varepsilon}\lim_{N\to+\infty}  \text{tr}(\log A_{-N,-b}(\varepsilon) -\log A_{-N,-b-1}(\varepsilon))=w(\varepsilon)
\end{equation}
and
\begin{equation}\label{lift1}
\det B_{a,p-b}(\varepsilon)=(-1)^nw(\varepsilon)\lim_{N,M\to +\infty} \frac{\det B_{-M,-b}(\varepsilon)\det B_{a,N}(\varepsilon)}{\det B_{-M,N}(\varepsilon)} \text{~~modulo~~} p,~  I^{\lceil Cp\rceil },~ \varepsilon^{\lceil Cp\rceil }.
\end{equation}

 \ 
 
 Notice that the l.h.s. of \eqref{w} \emph{does not depend} on $b\in \Z$. Indeed the shift $b\rightsquigarrow b+1$ of parameter $b$ is equivalent to the shift of variable $\varepsilon\rightsquigarrow \varepsilon+1$ in the ring $\widetilde{R}:=R[\varepsilon,\frac{1}{\varepsilon},\frac{1}{\varepsilon\pm 1},\frac{1}{\varepsilon\pm 2},...]$ and $T_\varepsilon$ is right-invariant with respect to the shift.
 
{\bf Proof.} It follows from Lemma 2 that $I$-adic limit 
\begin{equation}\label{rdet}
Q(\varepsilon)=\lim_{N,M\to +\infty} \frac{\det B_{-M,N}(\varepsilon)}{\det B_{-M,-b}(\varepsilon)\det B_{a,N}(\varepsilon)}
\end{equation}
exists.  It is clear that $\det B_{-M,-b}(\varepsilon),~\det B_{a,N}(\varepsilon)$ are invertible elements in $\hat{R}_I\hat{\otimes}_{\Z}\Q[[\varepsilon]]$ if $a,b>\max(|l_1|,|l_n|)$. Therefore, the  Weierstrass 
polynomial of $\frac{\det B_{-M,N}(\varepsilon)}{\det B_{-M,-b}(\varepsilon)\det B_{a,N}(\varepsilon)}$ is equal to $\varepsilon^n-w_{1,-M,N}\varepsilon^{n-1}+...+(-1)^{n}w_{n,-M,N}$. Let 
$w(\varepsilon)=\varepsilon^n-w_{1}\varepsilon^{n-1}+...+(-1)^{n}w_{n}$ be the  Weierstrass polynomial of $Q(\varepsilon)$. It is clear that $w(\varepsilon)=\lim_{N,M\to+\infty}(\varepsilon^n-w_{1,-M,N}\varepsilon^{n-1}+...+(-1)^{n}w_{n,-M,N})$.

The r.h.s. of (\ref{liftgen}) is an element in $\hat{R}_I\hat{\otimes}_{\Z}\Q[[\varepsilon]]$ so $w(\varepsilon)$ divides \begin{equation}
  \varepsilon^n\exp\left( T_{\varepsilon}\lim_{N\to+\infty}\text{tr}(\log A_{-N,1}(\varepsilon) -\log A_{-N,0}(\varepsilon)) \right) \label{limit}
\end{equation} 
\emph{ up to multiplication by an expression of the type $\exp(\phi)$ where $\phi\in pI+\epsilon^{\lceil Cp\rceil }I\hat{\otimes}\Q +I^{\lceil Cp\rceil }\hat{\otimes}\Q\subset I\hat{\otimes}\Q $ for all primes $p\gg 1$}. Using our assumption \eqref{ring_assumption} we conclude  that  $w(\varepsilon)$ divides \eqref{limit} rationally.

On the other hand,  we have $$
\varepsilon^n\exp\left( T_{\varepsilon}\lim_{N\to+\infty}\text{tr}(\log A_{-N,1}(\varepsilon) -\log A_{-N,0}(\varepsilon)) \right)=\varepsilon^n+ \text{elements of lower degree in }\varepsilon$$ which proves (\ref{w}).

The formula (\ref{lift1}) follows from (\ref{liftgen}) and (\ref{w}). $\square$

{\bf Remark 1.} The formula (\ref{rdet}) suggests to consider the Weierstrass polynomial $w(\varepsilon)$ as a regularized determinant of the infinite matrix $B(\varepsilon)$. One can show using methods of this and previous Sections that if $w(\varepsilon)\ne 0$, then $B(\varepsilon)$ is invertible 
but $B(\varepsilon)^{-1}$ belongs to a wider class of infinite matrices. In particular, $B(\varepsilon)^{-1}$ does not necessarily have only finitely many nonzero diagonals above the main diagonal, instead we have $\lim B(\varepsilon)^{-1}_{i,j}=0$ for $j-i\to +\infty$.

\section{p-Determinants of differential operators and the universal series $L(D)$}

Here and in the sequel we assume that $$R=\Z[t_{i,j}]$$ where $i,j=0,1,2,...$ and $I\subset R$ is an ideal generated by $t_{i,j}$ for $0\leq i\leq j$. Consider a differential operator
\begin{equation}\label{genD}
D=(x\partial_x-l_1)...(x\partial_x-l_n)+\sum_{\substack{0\leq j\leq m,\\ 0\leq i}}t_{i,j}x^i\partial_x^j
\end{equation}
for some $m\in\N,~l_1,...,l_n\in\Z$ and $l_1\leq l_2\leq...\leq l_n$. 

We denote by $D^{(p)}$ the reduction of $D$ modulo a prime number $p$. Let $V=\text{Span}_{R}(1,x,x^2,...)$ and $V^{(p)}=\text{Span}_{R}(x^p,x^{p+1},...)$. It is clear that the differential  operator $D^{(p)}$ acts on the completion of  $V$ with respect to the topology of Laurent series in $x$ and the subspace $V^{(p)}$ is invariant with respect to $D^{(p)}$. Therefore $D^{(p)}$ acts on the quotient space $V/V^{(p)}$ with  the  basis $1,x,x^2,...,x^{p-1}$. Let $Det_p(D)=\det D^{(p)} |_{V/V^{(p)}}$. It is clear that $Det_p(D)\in R/pR$. 

{\bf Theorem 3.} There exists an element $L(D)\in \hat{R}_I\hat{\otimes}_{\Z}\Q$ independent of $p$ such that $$Det_p(D)=L(D)~~~~ \text{mod}~ p~ \text{and}~ I^{\lceil Cp\rceil}$$ where $C>0$ is independent of $p$. 

{\bf Proof.}  Let $V_{\varepsilon}=\text{Span}_{R}(x^{i+\varepsilon};~i\in\Z)$. It is clear that $D$ acts on the  Laurent completion of $V_{\varepsilon}$. Let $B(\varepsilon)=(b_{ij}(\varepsilon)),~i,j\in\Z$ be the matrix of this operator with respect to the basis $\{x^{i+\varepsilon};~i\in\Z\}$. We have 
$$Dx^{j+\varepsilon}=\sum_{i\in\Z}b_{ij}(\varepsilon)x^{i+\varepsilon}$$
where
$$b_{ii}(\varepsilon)=(i+\varepsilon-l_1)...(i+\varepsilon-l_n)+\sum_{0\leq k\leq m}(i+\varepsilon)(i+\varepsilon-1)...(i+\varepsilon-k+1)t_{kk},$$
$$b_{ij}(\varepsilon)=\sum_{\max(0,j-i)\leq k\leq m}(j+\varepsilon)(j+\varepsilon-1)...(j+\varepsilon-k+1)t_{k+i-j,k}$$
for $i\ne j$.

It is clear that $B(\varepsilon)$ satisfies the properties listed in the beginning of the previous Section. Moreover, we have 
$Det_p(D)=\det B_{0,p-1}(0)$ and we can apply the formula (\ref{lift1}). $\square$

{\bf Theorem 4.} Let $$w_{-M,N}(\varepsilon)=\varepsilon^n-w_{1,-M,N}\varepsilon^{n-1}+...+(-1)^{n}w_{n,-M,N}$$ be the Weierstrass polynomial of $\det B_{-M,N}(\varepsilon)$ and 
$$w(\varepsilon)=\lim_{M,N\to\infty} w_{-M,N}(\varepsilon)$$ be its $I$-adic limit (see Theorem 2). Then we have 
\begin{equation}\label{L}
L(D)=(-1)^nw(0).
\end{equation}

{\bf Proof.} The formula (\ref{lift1}) gives
$$L(D)=\det B_{0,p-1}(0)=(-1)^nw(0)\lim_{N,M\to +\infty} \frac{\det B_{-M,-1}(0)\det B_{0,N}(0)}{\det B_{-M,N}(0)}.$$
Notice that $B_{-M,N}(0)$ is the  matrix of  the  differential operator $D$ restricted to the vector space spanned by $\{x^{-M},...,x^{-1},1,x,...,x^N\}$. This matrix is $2\times 2$ block triangular with 
matrices $B_{-M,-1}(0)$ and $B_{0,N}(0)$ on diagonal because the vector space spanned by $\{1,x,x^2,...\}$ is invariant with respect to $D$. Therefore $\det B_{-M,N}(0)=\det B_{-M,-1}(0)\det B_{0,N}(0)$ and we obtain (\ref{L}). $\square$

Let $\overline{R}$ be the smallest extension of $\hat{R}_I\hat{\otimes}_{\Z}\Q$ where the  Weierstrass polynomial $w(\varepsilon)$ completely factorizes\footnote{More precisely, $\overline{R}=\hat{R}_I\hat{\otimes}_{\Z}\Q[\varepsilon_1,...,\varepsilon_n]/\{\text{coefficients of}~\varepsilon^j,~j=0,...,n-1~~\text{in}~~ w(\varepsilon)-\prod_{j=1}^n(\varepsilon-\varepsilon_j)\}$.}. Let $\varepsilon_1,...,\varepsilon_n$ be all the  roots of $w(\varepsilon)$ in the ring $\overline{R}$. 

Notice that the formula (\ref{L}) can be written in the form
$$L(D)=\varepsilon_1...\varepsilon_n.$$

{\bf Lemma 6.} The differential equation $D g(x) =0$ has $n$ linearly independent formal solutions $g_1(x),...,g_n(x)$  satisfying the properties:

{\bf 1.}  $g_j(x)=x^{\varepsilon_j}\sum_{l\in\Z}r_{j,l}x^l$ for $j=1,...,n$ where $r_{j,l}\in\overline{R}$.

{\bf 2.}  $r_{j,l} \in I^{\max(0,\lceil \frac{C-l}{n}\rceil)}$ for  some constant $C$.                           

Moreover, if a formal non-zero solutions of the equation $D g(x) =0$ has a form $g(x)=x^{\delta}\sum_{l\in\Z}r_{l}x^l$ where $\delta=0$ mod $I$ and 
 there exist constants $C_1,C_2,C_3>0$ such that $r_{l}\in I^{\lceil C_2-C_1l\rceil}$ for $l<-C_3$, then $\delta=\varepsilon_j$ for some $j=1,...,n$.\footnote{One can show that in this case $g(x)$ is proportional to the solution $g_j(x)$ but we do not need this for our purposes.}

{\bf Proof.} We will approximate the infinite matrix $B(\varepsilon)$ of $D$ by finite matrices $B_{-M,N}(\varepsilon)$. Our solutions $g_1(x),...,g_n(x)$ are elements of the kernel of $B(\varepsilon)$ for $\varepsilon=\varepsilon_1,...,\varepsilon_n$ and we will obtain them as a limit of elements of the kernel of $B_{-M,N}(\varepsilon)$ for certain values of $\varepsilon$, when $M,N \to +\infty$. 

Recall that $w_{-M,N}(\varepsilon)=\varepsilon^n-w_{1,-M,N}\varepsilon^{n-1}+...+(-1)^{n}w_{n,-M,N}$ is the Weierstrass polynomial of $\det B_{-M,N}(\varepsilon)$. Let $\varepsilon_{1,-M,N},...,\varepsilon_{n,-M,N}$ be roots of $w_{-M,N}(\varepsilon)$. We have $\det B_{-M,N}(\varepsilon_{j,-M,N})=0$ for $j=1,...,n$. Let $g_{j,-M,N}(x)=x^{\varepsilon_{j,-M,N}}\sum_{-M\leq l\leq N}r_{j,l,-M,N}x^l$ be a non-zero element of 
the kernel of $B_{-M,N}(\varepsilon_{j,-M,N})$. In the $I$-adic limit $M,N\to \infty$ we get $\varepsilon_{j,-M,N}\to \varepsilon_j$ and $g_{j,-M,N}(x)\to g_j(x)$ where\footnote{One can prove the existence of this limit in the same way as in the proof of Lemma 2.} $g_j(x)$ generates the kernel of  $B(\varepsilon_j)$ and satisfies the property 1 by construction. The property 2  follows from the expression of a solution of  a  homogeneous system of linear equations in terms of minors of its matrix and the properties of the matrix $B(\varepsilon)$ listed in the beginning of Section 3. 

On the other hand, let $g(x)$ be a formal solution of the equation $D g(x) =0$ of the form $g(x)=x^{\delta}\sum_{l\in\Z}r_{l}x^l$ where $\delta=0$ mod $I$ and $r_{l}\in I^{\lceil C_2-C_1l\rceil}$ for $l<-C_3$. Choose an integer $L>\frac{1}{C_1}$. Introduce a new ring $R^{\prime}=\{\sum_{l\in\Z}a_lt^l,~a_l\in R~\text{for}~l>0,~a_l\in I^{-\lceil \frac{l}{L}\rceil}~\text{for}~l\leq 0\}$ where $t$ is a formal parameter. Define an ideal $I^{\prime}$ in $R^{\prime}$ by $I^{\prime}=tR^{\prime}$. We have $(I^{\prime})^L\cap R=I$. Let $D_0$ be a diagonal matrix for the operator $x^{\delta+l}\mapsto t^lx^{\delta+l},~l\in\Z$. Define $\tilde{B}(\delta)=D_0B(\delta)D_0^{-1}$. Notice that matrix elements of $\tilde{B}(\delta)$ belong to $R^{\prime}$ and off-diagonal elements belong to $I^{\prime}$. We have $\tilde{B}(\delta)\tilde{g}(x)=0$ where $\tilde{g}(x)=\sum_{l\in\Z}r_lt^{l+C_2}x^l$. Notice that  the  coefficients of $\tilde{g}(x)$ tend to zero in $I^{\prime}$-adic topology for both $l\to +\infty$ and $l\to-\infty$. 

If we reduce the equation  $\tilde{B}(\delta)\tilde{g}(x)=0$ modulo $(I^{\prime})^K$ and let $K\to+\infty$, we obtain a sequence of equations ~~~ $\tilde{B}_{-M,N}(\delta)\tilde{g}_{-M,N}(x)=0$~~~ with~~~ $M,N\to +\infty$ where ~~~$\tilde{g}_{-M,N}(x)=\sum_{-M\leq l\leq N}r_lt^{l+C_2}x^l$. Therefore, 
$\det \tilde{B}_{-M,N}(\delta)=\det B_{-M,N}(\delta)=0$ and $w_{-M,N}(\delta)=0$ modulo $(I^{\prime})^K$. In the limit $K\to\infty$ we get $w(\delta)=0$ and $\delta$ is equal to one of the roots of our Weierstrass polynomial\footnote{This can also be proved by constructing an inverse matrix $B(\delta)^{-1}$ in the case $w(\delta)\ne 0$, see Remark 1 in the end of Section 3.}.                      $\square$ 

{\bf Remark 2.}   If we were in the analytic framework (rather then in the algebraic one where we  actually are), we would say that  $e^{2\pi i\varepsilon_1},...,e^{2\pi i\varepsilon_n}$ are eigenvalues of  the  monodromy matrix $M_D$ of $D$ with respect to a small fixed loop around zero. We can write 
$$L(D)=\varepsilon_1...\varepsilon_n=\det \frac{\log M_D}{2\pi i}.$$
Therefore, we can reformulate the Theorems 3, 4  informally as
\begin{equation}\label{pdetinf}
Det_p(D)=\det \frac{\log M_D}{2\pi i}~ \text{~mod~}~ p \text{~~and~~} I^{\lceil Cp\rceil }.
\end{equation}
where $M_D$ is a "monodromy matrix" of $D$.

\section{Relation to monodromy matrix}

In this Section we provide an interpretation of the universal series $L(D)$ in the analytic framework, and make sense of the formula (\ref{pdetinf}).

Suppose that the parameter $m$ in the family  \eqref{family} satisfies inequality  $m\le n$, which means that our perturbation preserves the order of the differential operator $D$.

Let $R_{\text{an}}$ be the ring ${\mathcal O}_{z,z_0}$ of germs of analytic functions on a reduced irreducible complex analytic space $Z$ at some point $z_0\in Z$, and $I_{z_0}\subset R_{\text{an}}$ be the maximal ideal of $z_0$. Fix a homomorphism $$\mu:~ R \to R_{\text{an}}$$ such that $\mu(I)\subset I_{z_0}$. This homomorphism maps $t_{ij}\in R$ to germs of analytic functions at $z_0$ which we denote by $t_{ij}(z)$. We have $t_{ij}(z_0)=0$ for $0\leq i\leq j$. 

Assume that there exists $\rho >0$ such that for all $j=0,...,m$ the series $\sum_{i=0}^{\infty} t_{ij}(z) x^i$ is absolutely convergent for all $x\in \C,~~~|x|<\rho$ and all $z$ sufficiently close to $z_0$. Under our assumptions,  the formula (\ref{genD}) gives a holomorphic family $D_z$ of differential operators\footnote{More formally, we have $D_z=\mu(D)$.} of order $n$ in the open disc $|x|<\rho$.  The symbol of $D_{z_0}$ vanishes only at $x=0$.  Pick a number $\rho^{\prime}$ such that $0<\rho^{\prime}<\rho$. Then for $z$ sufficiently close to $z_0$ the symbol of $D_z$ does not vanish in the annulus $A_{\rho^{\prime},\rho}=\{x;~\rho^{\prime}\leq |x|\leq \rho\}$. 

Solutions of the equation $D_z g(x)=0$ form a local system of rank $n$ on $A_{\rho^{\prime},\rho}$. Denote by $M_{D_z}$ the monodromy matrix (which is an element of $\text{GL}(n,\C)$ defined up to conjugation) of this system along the loop in $A_{\rho^{\prime},\rho}$ (e.g. $x(\theta)=\frac{\rho^{\prime}+\rho}{2} e^{i\theta},~\theta\in [0,2\pi]$). Notice that $M_{D_z}$ depends holomorphically on $z$, and is a unipotent operator for $z=z_0$. Therefore, we have a canonical choice of $\log M_{D_z}$ for $z$ close to $z_0$. 

We make further assumptions that

{\bf 1.} For an open dense subset $U\subset Z$ the spectrum of $\log M_{D_z}$ is simple,

{\bf 2.} Eigenvalues of $\log M_{D_z}$ are univalued functions on $U$.

The first assumption can be achieved by an embedding of $Z$ into a larger germ. The second assumption is achieved by passing to an appropriate ramified cover. 

These two assumptions imply that the equation $D_zg(x)=0$ has $n$ linearly independent  solutions in $A_{\rho^{\prime},\rho}$ of the form 
$$g_j^{\text{an}}(x)=x^{\varepsilon_j^{\text{an}}}\sum_{l\in\Z} r_{j,l}^{\text{an}} x^l$$
where $\varepsilon_j^{\text{an}}$, $j=1,...,n$ are pairwise distinct, $\varepsilon_j^{\text{an}}\in I_{z_0}$, $r_{j,l}^{\text{an}}\in R_{\text{an}}$, and $$g_j^{\text{Laur}}(x)=x^{-\varepsilon_j^{\text{an}}}g_j^{\text{an}}(x)=\sum_{l\in\Z} r_{j,l}^{\text{an}} x^l$$ is a Laurent series in $x$ which is uniformly bounded in $A_{\rho^{\prime},\rho}$ for all $z$ close to $z_0$.  

{\bf Theorem 5.} Under above assumptions, the homomorphism $\mu:~R\to R_{\text{an}}$ extends to a homomorphism $\overline{\mu}:~\overline{R}\to R_{\text{an}}$  and $\varepsilon_j^{\text{an}}=\overline{\mu}(\varepsilon_j)$, $j=1,...,n$ up to some permutation from $S_n$.\footnote{One can show that solutions of the equation $D_z g(x)=0$  are equal to $g_j^{\text{an}}=\overline{\mu}(g_j)=x^{\overline{\mu}(\varepsilon_j)}\sum_{l\in\Z}\overline{\mu}(r_{j,l})x^l,~~~j=1,...,n$, up to multiplication by a constant depending on $z$ and $j$, where $g_1,...,g_n$ are defined in Lemma 6.} 

The proof is based on the following

{\bf Lemma 7.} Let $R_{\text{for}}$ be a ring containing $\Q$, without zero divisors, and complete with respect to a maximal ideal $I_{\text{for}}\subset R_{\text{for}}$. Fix a homomorphism $\mu_{\text{for}}: \overline{R} \to R_{\text{for}}$ such that $\mu_{\text{for}}(I)\subset I_{\text{for}}$ and $\mu_{\text{for}}(\varepsilon_j)$ are pairwise distinct. Let $D_{\text{for}}=\mu_{\text{for}}(D)$.                        

If a non-zero formal solutions of the equation $D_{\text{for}} g(x) =0$ has a form $g(x)=x^{\delta}\sum_{l\in\Z}r_{l}x^l$ where $\delta=0$ mod $I_{\text{for}}$ and  $r_{l}\in I^{\lceil -Cl\rceil}$ for a constant $C>0$ and $l<0$, then $\delta=\mu_{\text{for}}(\varepsilon_j)$ for some $j=1,...,n$.

{\bf Proof.} Omitted because it is similar to the proof of the second part of Lemma 6.  $\square$

{\bf Proof of the Theorem.}  We assume (without loss of generality) that the germ of $Z$ at $z_0$ is embedded analytically in $\C^N$, so we can speak about distance in $Z$. There exist constants $0<C_1<1,~0<C_2$ such that for sufficiently small $\delta=|z-z_0|$ operator $D_z$ has  a  non-vanishing symbol in the annulus $\delta^{C_1}\leq |x|\leq \rho$ and the holonomy matrix for the shortest path connecting $\frac{\rho}{2}$ and any point of the circle 
$\{e^{i\theta}\delta^{C_1};~0\leq\theta\leq 2\pi\}$ is bounded by $\delta^{-C_2}$. 

Let a multivalued solution of the equation $D_zg(x)=0$  has the form $g_j^{\text{an}}(x)=x^{\varepsilon_j^{\text{an}}}g_j^{\text{Laur}}(x)$ where
$g_j^{\text{Laur}}(x)=\sum_{l\in\Z} r_{j,l}^{\text{an}} x^l.$  The uniform bound on $g_j^{\text{Laur}}(x)$ in $A_{\rho^{\prime},\rho}$ implies 
that $g_j^{\text{Laur}}(x)$ extends analytically to the annulus $A_{\delta^{C_1},\rho}$, and satisfies the bound 
 $$|g_j^{\text{Laur}}(e^{i\theta}\delta^{C_1})|<\text{const}~\delta^{-C_2}.$$ We can write
$$r_{j,l}^{\text{an}}=\frac{1}{2\pi i} \oint_{|x|=\delta^{C_1}} g_j^{\text{Laur}}(x)  x^{-l} \frac{dx}{x}=\frac{1}{2\pi i} \oint_{|x|=\rho} g_j^{\text{Laur}}(x)  x^{-l} \frac{dx}{x}.$$
This implies the bound
$$|r_{j,l}^{\text{an}}|\leq \begin{cases}~ \text{const}~\delta^{-C_2-C_1l},~~~\text{if}~~~l<0, \\ ~
 \text{const}~\rho^{-l},~~~\text{if}~~~l\geq 0.\end{cases}$$
This implies $r_{j,l}^{\text{an}}\in  I^{\lceil -Cl\rceil}$ for some $C>0$ and $l<0$. Using Lemma 7 we conclude that $\varepsilon_l^{\text{an}}=\mu(\varepsilon_l)$ and the Theorem follows.   $\square$

{\bf Corollary.} The formal completion of the germ of analytic function 
$$z \mapsto \det \frac{\log M_{D_z}}{2\pi i}$$
coincides with $L(D_z)$.

{\bf Remark 3.} In the case $m>n$ the analytic interpretation of the series $L(D)$ is not clear. For example, consider the differential operator $(x\partial_x-l_1)...(x\partial_x-l_n)+t\partial_x^m$ with small parameter $t$ and $m>n$. One can show that in this case  $\varepsilon_j=0$ for all $j$.

\section{Examples}

{\bf Example 1.} Let $$D=\frac{d}{dx}x(x-1)(x-t)\frac{d}{dx}+x=x(x-1)(x-t)\frac{d^2}{dx^2}+(3x^2-2(t+1)x+t)\frac{d}{dx}+x.$$
This operator can be written as $D=-(t+1)D_0$ where 
$$D_0=x\frac{d}{dx}\Big(x\frac{d}{dx}+1\Big)-\frac{1}{t+1}\Big(x^3\frac{d^2}{dx^2}+3x^2\frac{d}{dx}+x+tx\frac{d^2}{dx^2}+t\frac{d}{dx}\Big)$$
It is clear that $Det_p(D)=(-1)^p(1+t)^pDet_p(D_0)$ as determinants of $p\times p$ matrices. We have 
$(-1)^p(1+t)^p=-1$ mod $p$ and $t^p$ for $p>2$. The differential operator $D_0$ has a form (\ref{genD}). Let 
$e^{\pm \pi i\lambda}$ be eigenvalues of the monodromy matrix $M_{D_0}$. We have 
$$Det_p(D)=(-1)^p(1+t)^pDet_p(D_0)=-Det_p(D_0)=\lambda^2 \text{~~mod~~}p,~t^{p-1}$$
which partially proves\footnote{Our methods give the proof  of this identity modulo $p,~t^{\lceil \frac{p}{2}\rceil}$.} the Claim from the Introduction. Notice that $\lambda$ is a solution of the equation\footnote{This is a well known fact in the theory of Heun functions. See for example  {\url{https://dlmf.nist.gov/31}} and references therein.}
\begin{equation}\label{conteq}
1-\frac{\frac{t}{(t+1)^2} \frac{(\lambda+2)^2}{(\lambda+1)(\lambda+3)}}{1-\frac{\frac{t}{(t+1)^2} \frac{(\lambda+3)^2}{(\lambda+2)(\lambda+4)}}{1-\frac{\frac{t}{(t+1)^2} \frac{(\lambda+4)^2}{(\lambda+3)(\lambda+5)}}{1-\dots}}}~+~ 1-\frac{\frac{t}{(t+1)^2} \frac{(\lambda+1)^2}{\lambda(\lambda+2)}}{1-\frac{\frac{t}{(t+1)^2} \frac{\lambda^2}{(\lambda-1)(\lambda+1)}}{1-\frac{\frac{t}{(t+1)^2} \frac{(\lambda-1)^2}{(\lambda-2)(\lambda)}}{1-\dots}}}=1.
\end{equation}
The power series expansion (\ref{lambda}) can be found from this equation.  

As we already mentioned in the Introduction, the coefficients of both power series $\lambda (t)$ and $h(t)=\lambda (t)^2$ have very large denominators. For example, the coefficient of $h(t)$ at $t^{20}$ is 

$$\frac{8488029376721954239470359915918389605555923159104304646275144404309630883}{2562649082054740025416814617780707052281611121021293690880000000000000000}$$
with the  denominator equal  to $2^{37}\cdot 3^{24}\cdot 5^{16}\cdot 7^{14}\cdot 11^{10}\cdot 13^8\cdot  17^4\cdot 19^2$. Experiments suggest that the  $n$-th coefficient of $h(t)$ for $n\geq 2$ has a form
$$\text{sqn}_n\cdot \frac{b_n}{\prod_{p\leq n}p^{\alpha_p(n)}}$$
where $p$ are prime numbers, $b_n\geq 1$ and $(b_n,p)=1$ for $p\leq n$. 

Moreover, $b_n=e^{\frac{1}{2}n^2(1+o(1))}$ and\footnote{The following formulas has exceptions in the case $\alpha_3(n)$,  $n=6,7,8$.}

$$\alpha_2(n)=\text{ord}_2(2^{n-1}n!),$$
$$\alpha_p(n)=\max_{k\geq 1}(k(n+1-p^k)),~~~3\leq p\leq n,$$
$$\text{sqn}_n=(-1)^n,~~~n\geq 6.$$

This implies that $h(t)$ does not satisfy a non-trivial Picard-Fuchs equation (which also follows from the monodromy formula). It was speculated in 
Section 4.3 of \cite{K} that a certain algebra of $p$-determinants (denoted by $P_k\subset k_{\infty}$) is contained in the algebra of periods. The 
present study of $p$-determinants was motivated by an attempt to verify this conjecture from loc.cit. in the parametric case. Now we see that it is disproved and should be replaced by something new in light of our main result relating $p$-determinants and monodromy.

{\bf Example 2.} Let $$D_1=x(x-1)(x-t)\frac{d^2}{dx^2}+\frac{1}{2}\Big(x(x-1)+x(x-t)+(x-1)(x-t)\Big)\frac{d}{dx}+1.$$
This differential operator can be written as
$$D_1=\Bigg(\sqrt{x(x-1)(x-t)}\frac{d}{dx}\Bigg)^2+1$$
and therefore the differential equation $D_1f(x)=0$ has solutions $$f_{\pm}(x)=\exp\Bigg(\pm i\int_{x_0}^x\frac{dx}{\sqrt{x(x-1)(x-t)}}\Bigg).$$
A monodromy matrix of $D_1$ with respect to an arbitrary loop has eigenvalues of the form $\exp(\pm 2\pi i\lambda(t))$ where 
$$\lambda(t) =\frac{1}{2\pi }\int_{\Gamma}\frac{dx}{\sqrt{x(x-1)(x-t)}} \text{~~~modulo~~~} \Z$$
and $\Gamma$ is a cycle on the elliptic curve $y^2=x(x-1)(x-t)$. In particular, for a loop around $0,t$ if $|t|<1$ we have  
$$ \lambda(t)=\sum_{k=0}^{\infty}\frac{(2k)!^2}{2^{4k}k!^4}t^k=1+\frac{1}{4}t+\frac{9}{64}t^2+\frac{25}{256}t^3+\frac{1225}{16384}t^4+\frac{3969}{65536}t^5+...$$

On the other hand we can write $D_1$ in the form 
$$D_1=-\Big(x\frac{d}{dx}+1\Big)\Big(x\frac{d}{dx}-1\Big)+x^3\frac{d^2}{dx^2}+\frac{3}{2}x^2\frac{d}{dx}-t(x^2-x)\frac{d^2}{dx^2}-t(x-\frac{1}{2})\frac{d}{dx}.$$
Therefore, we have 
$$Det_p(D_1)=\Big(\frac{1}{4}t+\frac{9}{64}t^2+\frac{25}{256}t^3+\frac{1225}{16384}t^4+\frac{3969}{65536}t^5+...\Big)^2=\frac{1}{16}t^2+\frac{9}{128}t^3+\frac{281}{4096}t^4+...$$
modulo $p,~t^{\lceil Cp\rceil}$. Notice that denominators of coefficients are powers of $2$ in this case.

We see that in this example the $p$-determinant is related to the series $\lambda(t)$ satisfying a Picard-Fuchs equation. On the other hand,  the same series $\lambda(t)$ appeared in \cite{M} in relation to the Hasse invariant of the elliptic curve $y^2=x(x-1)(x-t)$. Notice that its Hasse invariant is the same as for the twisted curve $E=\{(x,y);~~y^2=-x(x-1)(x-t)\}$. Moreover, it is essentially equivalent to the $p$-curvature of the 
connection on the trivial bundle of rank one  on $E$ given by $d+\frac{dx}{y}$. The differential operator $D_1$ is the direct image of this connection with respect to the map $E\to \P^1$ given by $(x,y)\mapsto x$ on $\A^1\subset \P^1$. It was explained in \cite{K} that characteristic polynomial of a $p$-curvature of a differential operator on $\A^1$ can be calculated in terms of $p$-determinants. This gives an alternative explanation to the  Manin's original observation.

{\bf Remark 4.} One can show that the arithmetic support (see \cite{K}, Section 3.2) of differential operators in both Examples 1, 2 is the zero locus of
$$\tilde{x}^p(\tilde{x}^p-1)(\tilde{x}^p-t^p)\tilde{y}^{2p}+Det_p(D)$$
considered as an element of $ \Z/p\Z[t][\tilde{x}^p,\tilde{y}^p]$,  see notations in  loc.cit.

\addcontentsline{toc}{section}{Acknowledgement}

\section*{Acknowledgement}

We are grateful to Don Zagier for useful discussions in the early stages of this project. A.O. is grateful to IHES for invitations and excellent working atmosphere.


\addcontentsline{toc}{section}{References}




\begin{thebibliography}{99}

\bibitem{K} Maxim  Kontsevich, Holonomic D-modules and positive characteristic. Japanese Journal of Mathematics, March 2009, Volume 4, Issue 1, pp 1-25

\bibitem{M} Yu. I. Manin,  The Hasse-Witt matrix of an algebraic curve.  AMS  Translations,  ser.  2,  vol.  45  (1965), pp 245-264




\end{thebibliography}
\end{document}